\documentclass[12pt]{article}
\usepackage{amsfonts, amsmath, amsthm, amssymb}

\newtheorem{theorem}[equation]{Theorem}
\newtheorem{lemma}[equation]{Lemma}
\newtheorem{prop}[equation]{Proposition}
\newtheorem{question}[equation]{Question}

\def\CC{\mathbb{C}}
\def\FF{\mathbb{F}}
\def\PP{\mathbb{P}}
\def\QQ{\mathbb{Q}}
\def\ZZ{\mathbb{Z}}

\DeclareMathOperator{\Cl}{Cl}
\DeclareMathOperator{\Gal}{Gal}
\DeclareMathOperator{\ord}{ord}
\DeclareMathOperator{\Res}{Res}

\newcounter{fixmectr}

\begin{document}

\title{Quantum computation of zeta functions of curves}
\author{Kiran S. Kedlaya \\ Department of Mathematics, Room 2-165 
\\ Massachusetts
Institute of Technology \\ 77 Massachusetts Avenue \\
Cambridge, MA 02139 \\
\texttt{kedlaya@math.mit.edu}}
\date{November 29, 2005}

\maketitle

\begin{abstract}
We exhibit a quantum algorithm for determining the zeta function
of a genus $g$ curve over a finite field $\FF_q$, which is 
polynomial in $g$ and $\log(q)$.
This amounts to giving an algorithm to produce provably random
elements of the class group of a curve, plus a recipe for
recovering a Weil polynomial from enough of its cyclic
resultants. The latter effectivizes a result of Fried in a restricted setting.
\end{abstract}

\section{Introduction}

Given a curve $C$ (assumed to be smooth, projective and geometrically
irreducible) over a finite field $\FF_q$ with $q = p^a$ for some prime $p$,
the zeta function of $C$ has the form
\[
Z(C,t) = \exp \left( \sum_{n=1}^\infty \frac{T^n}{n} \#C(\FF_{q^n}) 
\right) = \frac{P(t)}{(1-t)(1-qt)}
\]
for some polynomial $P(t) \in \ZZ[t]$ of degree $2g$ with $P(0) = 1$.
The determination of $P(t)$ is an active problem in algorithmic
number theory, in part because of practical connections to cryptography
(especially when $C$ is an elliptic curve, or more generally
a hyperelliptic curve).
For $g$ fixed, the approach introduced by Schoof \cite{schoof} (compute
$P(t)$ modulo many small primes) gives an algorithm
which is polynomial
in $\log(q)$ but exponential in $g$, as shown by Pila \cite{pila}
and Adleman-Huang \cite{adleman-huang}.
(A streamlined form of Schoof's algorithm, incorporating improvements due
to Atkin, Elkies, et al., turns out to be usable in practice
for $g=1$ and perhaps for $g=2$, but for larger $g$ the
algorithm is highly impractical.)
On the other hand, 
imitating Dwork's proof of the rationality of zeta functions \cite{dwork}
yields an algorithm which is polynomial in $p$, $g$ and $\log_p(q)$,
as observed by Lauder and Wan \cite{lauder-wan}. (The latter is 
also not practical, but related ``cohomological'' techniques have
proven more tractable; see \cite{denef-vercauteren} for the current state
of the art.)

However, a single algorithm for computing $P(t)$ in time polynomial 
\emph{both} in $g$ and $\log(q)$ remains elusive. Thus any sign that this
problem might be ``easy'' has some relevance; the main result of this note
(originally written as an addendum to \cite{kedlaya})
is one such sign, if only an indirect one.
\begin{theorem}
There is a quantum algorithm for computing the numerator
$P(t)$ of the zeta function,
which is polynomial time in $g, \log(q)$.
(See Section~\ref{sec:prob} for conventions regarding probabilistic
algorithms.)
\end{theorem}
Implicit in the statement of the theorem is the choice of
a mechanism for inputting arbitrary curves, such that the length of the input
is polynomial in the genus. We will be more explicit about the choice we have
in mind in Section~\ref{sec:Jacobians}; 
however, if the reader prefers to substitute
a polynomial time equivalent alternate choice,
this will of course not affect the truth of the theorem.

The components of the algorithm specified in Theorem~1 will be
described in the subsequent sections of the paper. It may be worth
pointing out here some components that may have some interest on their own:
a method for producing generators of the Jacobian group of a curve
over a finite field with provably high probability (Lemma~\ref{L:generators}),
and a method for recovering a Weil polynomial from a few of its cyclic
resultants (Section~\ref{sec:zeta}).

\section{Conventions for probabilistic algorithms}
\label{sec:prob}

Before proceeding, it will be helpful to fix some conventions about
probabilistic algorithms.

Given a real number $b \in (0,1)$,
we define a \emph{Las Vegas algorithm} to be an algorithm that, given a stream of outputs
of a ``fair coin'' (a/k/a a Bernoulli trial with probability 1/2), accomplishes its specified
goal with probability at least $1-b$ and reports failure with probability $b$.
As long as $b$ is fixed, its exact value is not critical, as 
the success probability of a Las Vegas algorithm can be boosted simply by repeated
invocation. This analysis is standard (and easy), but it will be useful for us to record
it explicitly: in terms of the success probability $a = 1-b$, in case $a \leq 1/2$, then
after two invocations, the success probability is 
\[
1 - (1-a)^2 = 2a - a^2 = a(2 - a) \geq 3a/2.
\]
In particular, one can boost the success probability from $a$ to $1/2$ with at most
\[
2^{\lceil \log_{3/2} (2/a) \rceil}
\leq 2^{2 \log_2 (2/a) + 1}
= \frac{8}{a^2}
\]
invocations, and from there up to any fixed higher
value by multiplying the number of invocations by a suitable fixed factor.
For instance, to get to success probability $3/4$, it suffices to perform
$16/a^2$ invocations.
(By the same token, 
it is sometimes more convenient to use Bernoulli trials of different probabilities,
e.g., to sample uniformly from a finite set; one can simulate such trials with a fair
coin up to any fixed failure probability.)

Given a real number $b \in (1/2,1)$,
we define a \emph{Monte Carlo algorithm} to be an algorithm that, given a stream of outputs
of a fair coin, accomplishes its specified
goal with probability at least $1-b$ but may yield any outcome otherwise.
Because of the nature of quantum mechanics, all quantum algorithms must be
regarded as Monte Carlo algorithms.
Again, one can decrease the error probability $b$ below any fixed cutoff, this time by
performing a fixed number of invocations and retaining the answer returned
most often. This analysis is standard, and it will not be useful for us 
to record it explicitly,
so we omit it.

\section{Black box groups}

Our quantum algorithm for computing zeta functions reduces the problem
to the determination of the order of certain ``black box groups''.
Before proceeding to the specific groups in question (groups of
rational points on Jacobian varieties), we first recall a bit of the
formalism of black box groups and cite the result about them we will be using.
Note that this formalism makes sense within any of the standard computing paradigms (e.g., deterministic, Las Vegas, Monte Carlo, or quantum).

A \emph{black box group with unique encodings},
in the sense of
Babai and Szemer\'edi \cite{babai-szemeredi}, consists of an
$n$-element subset $T$ of $\{0,1\}^m$ for some $m$ and $n$, and an oracle which has the following
properties, for some (unknown) subset $S \subseteq \{0,1\}^m$ containing $T$, and some (unknown) bijective map $f: S \to G$ from $S$ to
a group $G$ generated by $f(T)$.
\begin{enumerate}
\item[(a)]
Given $x,y \in S$, the oracle can determine $z \in S$ such that
$f(z) = f(x) f(y)$ in $G$. 
\item[(b)]
Given $x \in S$, the oracle can determine $y \in S$ such that
$f(y) = f(x)^{-1}$ in $G$.
\end{enumerate}
We may also speak of this data as a ``black box presentation of $G$ with
unique encodings''; its input length for complexity purposes is taken to be $mn$. 
Compare this definition with that of a
``black box group'' without further qualification:
in that case $f$  is only required to be
surjective, and the oracle is required to be able to determine,
given $x \in S$, whether $f(x)$ is the identity element of $G$.

We are now ready to invoke the necessary input from the theory of 
quantum computing.
\begin{lemma} \label{L:quantum}
Given a Monte Carlo black box group presentation with unique encodings
$f: S \to G$ of an abelian (or even solvable) group $G$,
of input length $mn$, 
there is a quantum algorithm, running in time polynomial in $mn$,
for computing the order of $G$.
\end{lemma}
\begin{proof}
  See Watrous \cite{watrous1}, \cite{watrous2}; 
the technique extends Shor's application of Fourier transform
methods to the factoring and discrete logarithm problems
\cite{shor}.
\end{proof}

\section{Algebraic curves}

Since our intended reader is not necessarily an expert in algebraic geometry,
we include here a synopsis of some relevant facts. 
For a fuller treatment, see \cite{fulton} or
\cite[Chapter~IV]{hartshorne}.

By a \emph{curve} over a perfect field $k$, 
we will always mean a smooth, projective, geometrically
irreducible variety $C$ of dimension 1 over $k$. To each such curve we can associate the
field $K(C)$ of rational functions on $C$; this is a field of
transcendence degree 1 over $k$, in which $k$ is relatively
algebraically closed. In fact, the functor $C \mapsto K(C)$ is an equivalence
between curves and such fields.
Let $\overline{k}$ denote the algebraic closure of $k$,
and let $C(k)$ and $C(\overline{k})$ denote the sets of
$k$-rational and $\overline{k}$-rational points, respectively,
on $C$.

A \emph{divisor} on $C$ is a formal sum
\[
D = \sum_{P \in C(\overline{k})} c_P(P) \qquad (c_P \in \ZZ),
\]
invariant under the action of 
$\Gal(\overline{k}/k)$ induced by the Galois action
on $C(\overline{k})$,
in which $c_P = 0$ for all but finitely many $P$. That last condition means
that the sum $\sum_P c_P$ is well-defined; it is called the
\emph{degree} of $D$ and denoted $\deg(D)$.

We point out three special types of divisors.
We refer to the sum over a single Galois orbit on $C(\overline{k})$, with
all coefficients 1, as a \emph{prime divisor}; the group of divisors is
freely generated by the prime divisors.
For $f \in K(C)^*$ and $P \in C(\overline{k})$, let 
$\ord_P(f)$ denote the order of vanishing (positive, negative, or zero)
of $f$ at $P$.
Define the divisor $(f) = \sum_P \ord_P(f) (P)$; any divisor of this form
is called a \emph{principal divisor}.
Similarly, for $\omega$ a nonzero 1-form on $C$, we may define
$\ord_P(\omega)$ as the order of vanishing, and define the divisor
$(\omega) = \sum_P \ord_P(\omega) (P)$; any divisor of this form is 
called a \emph{canonical divisor}. Note that if $D$ is a principal
divisor, then $\deg(D) = 0$, whereas if $D$ is a canonical divisor,
then $\deg(D) = 2g-2$, where $g$ is the genus of $C$ (by the Riemann-Roch
theorem; see below).
We write $D_1 \sim D_2$ to mean that $D_1 - D_2$ is a principal divisor;
this is clearly an equivalence relation. Note that the ratio of two
1-forms is a rational function, so any two canonical divisors are equivalent.

A divisor $D = \sum_P c_P(P)$ is \emph{effective} if $c_P \geq 0$ for all $P$;
we write $D_1 \geq D_2$ to mean that $D_1-D_2$ is effective.
For $D$ effective, we necessarily have $\deg(D) \geq 0$ (but not
conversely).
For $D$ a divisor on $C$, let $L(D)$ be 
the set of functions $f \in K(C)$ such that
$(f) + D \geq 0$, together with the zero function.
The set $L(D)$ is a vector space over $k$; let
$\ell(D)$ be the dimension of that space. Note that
$\ell(D) = 0$ whenever $\deg(D) < 0$.
The main theorem governing $\ell(D)$ is the Riemann-Roch theorem,
whose statement is the following.
\begin{prop}[Riemann-Roch theorem]
For any divisor $D$ on $C$,
\[
\ell(D) = \deg(D) + 1 - g + \ell(K-D).
\]
\end{prop}

The \emph{class group} $\Cl(C)$ is defined as the group of divisors of degree
zero, modulo the subgroup of principal divisors; it can be identified with
the $k$-rational points of a certain $g$-dimensional abelian variety
$J$, the so-called \emph{Jacobian variety} of $C$. Over a finite field,
the order of 
$\Cl(C)$ is closely related to the zeta function, by the following formula
(for which see, e.g., \cite[Section 14]{milne}).
\begin{prop} \label{P:Jacobian}
Suppose $k = \FF_q$; let $C_n$ denote the base change of $C$ to
$\FF_{q^n}$.
Let $P(t)$ be the numerator of the zeta function of $C$. Then
$\deg(P) = 2g$, and if we factor
$P(t) = (1 - r_1 t)\cdots(1 - r_{2g} t)$ with $r_1, \dots, r_{2g} \in \CC$, then
\[
\#\Cl(C_n) = \prod_{i=1}^{2g} (1 - r_i^n).
\]
\end{prop}
For this reason, computing the order of $\Cl(C)$ when $k$ is finite is
key to our quantum algorithm for computing zeta functions. The order is further
controlled by the Riemann hypothesis for curves
(see \cite[Chapter~X]{lorenzini} for a not-too-technical
treatment).
\begin{prop} \label{P:Weil bound}
With notation as in Proposition~\ref{P:Jacobian}, 
$|r_i| = q^{1/2}$ for $i=1, \dots, 2g$. In particular,
\begin{gather*}
q^n - 2g q^{n/2} \leq \#C(\FF_{q^n}) \leq q^n + 2g q^{n/2} \\
q^{ng/2} (\sqrt{q} - 1)^{ng} \leq \#\Cl(C_n) \leq q^{ng/2} 
(\sqrt{q} + 1)^{ng}.
\end{gather*}
\end{prop}

We will exploit the Riemann hypothesis via the following lemma.
\begin{lemma} \label{L:count primes}
For $e$ a positive integer, the number of prime divisors of degree $e$
on $C$ is at least
\[
\frac{1}{e} (q^e (1-q^{-1}) - 4gq^{e/2}).
\]
\end{lemma}
\begin{proof}
It suffices to count elements of $C(\FF_{q^e})$, subtract elements of
$C(\FF_{q^i})$ for all proper divisors $i$ of $e$, and then
divide by $e$. By Proposition~\ref{P:Weil bound},
this count can be bounded below by
\[
\frac{1}{e}(q^e - 2gq^{e/2} - \sum_{i<g, i|g} (q^i + 2g q^{i/2})).
\]
If $e=1$, there is no sum at right, so we obtain
$q - 2g q^{1/2}$ as the lower bound, which implies the desired bound. If $e=2$, the bound is
\[
\frac{1}{2}(q^2 - 2gq - q - 2gq^{1/2})
\geq \frac{1}{2}(q^2 (1 - q^{-1}) - 4gq).
\]
Otherwise we may dominate
$\sum_{i<g,i|g} q^i$ by $\sum_{i=1}^{e-2} q^i \leq q^{e-1}$,
and we may dominate
$\sum_{i<g,i|g} 2g q^{i/2}$ by $\sum_{i=1}^{\lfloor e/2 \rfloor} 2g q^{i/2}
\leq 2g q^{e/2}$.
This yields the lower bound
\[
\frac{1}{e}(q^e - 2g q^{e/2} - q^{e-1} - 2g q^{e/2}) =
\frac{1}{e}(q^e (1-q^{-1}) - 4g q^{e/2}).
\]
\end{proof}

\section{Representing elements of class groups}
\label{sec:rep}

In the notation of the previous section, we collect here some observations
about representing elements of $\Cl(C)$. 

We first note that elements can be
represented in a compact form. Let $U$ be a divisor with
$\deg(U) = 1$. Given a divisor $D$ with $\deg(D) = 0$, we have by
Riemann-Roch
\[
\ell(D + mU) = m + 1 - g + \ell(K-D-mU) \geq m + 1 - g;
\]
in particular, if $m \geq g$, then $\ell(D + mU) > 0$,
so that $D + mU \sim E$ for some effective divisor $E$. In other
words, for any fixed $m \geq g$,
every element of $\Cl(C)$ can be represented as $E - mU$ for some
effective divisor $E$ of degree $m$.

The representations of elements of $\Cl(C)$ in the form $E - gU$,
for $E$ effective of degree $g$, are unique ``generically'' but not always;
since we will need to generate random elements of $\Cl(C)$, it will be
useful to have representations which are 
uniformly distributed across $\Cl(C)$. Namely, if $\deg(D) = 0$ and
$m \geq 2g-1$, 
we have $\deg(K - D - mU) < 0$, so Riemann-Roch yields
$\ell(D + mU) = m+ 1 - g$. In particular, if $k = \FF_q$,
then each element of $\Cl(C)$ is represented by exactly $q^{m+1-g}$ number
of divisors of the form $E - mU$, for $E$ effective of degree $m$.

Finally, we note that in case there exists a rational point 
$O \in C(k)$, we can represent
elements of $\Cl(C)$ in a canonical form. Namely, in this case, if
$\deg(D) = 0$, then
\begin{align*}
\ell(D + (m-1)(O)) &\leq \ell(D + m(O)) \\
&= m + 1 - g + \ell(K-D-m(O)) \\
&\leq m+1-g+\ell(K-D-(m-1)(O)) \\
&= \ell(D + (m-1)(O)) + 1.
\end{align*}
Hence if $m$ is the smallest nonnegative integer for which 
$\ell(D + m(O)) > 0$, then $m \leq g$ (as above) and $\ell(D + m(O)) = 1$.
In other words, for this choice of $m$ (which depends on $D$),
there is a \emph{unique} effective divisor $E$ 
with $D + m(O) \sim E$.

 \section{Computing in class groups}
\label{sec:Jacobians}

We now make some remarks about the protocols we have in mind for
inputting and computing on algebraic curves, starting with what
constraints on these protocols are imposed by the demands of our algorithm.
Note that we will make liberal use of factorization of monovariate polynomials
over finite fields, so our algorithms will be Las Vegas rather than
deterministic.

Let $C$ be a curve (which as usual
is smooth, projective, and geometrically irreducible)
of genus $g$ over $\FF_q$; for $n$ a positive integer, let $C_n$
be the base change of $C$ to $\FF_{q^n}$.
For the proof of Theorem~1 we will
need an algorithm to compute $\#\Cl(C)$ in time polynomial in
$g$ and $\log(q)$.
Using Lemma~\ref{L:quantum}, we see that it is enough
to exhibit a Monte Carlo black box presentation with unique encodings
of $\#\Cl(C)$, of input length
bounded by a polynomial in $g$ and $\log(q)$; in fact, our oracular operations
will be Las Vegas and not just Monte Carlo.
Beware that for technical reasons, we will eventually have to restrict to the
situation where $q$ is ``not too small'' compared to $g$; however, that
restriction will not be relevant in this section. (It will also be dropped out
in the course of proving Theorem~1.)

We now proceed to describing our input protocol and the construction
of the black box presentation of $\#\Cl(C)$, except for producing a
generating set; we defer that construction to the next section.
To begin with, we will input $C$ by specifying a homogeneous polynomial in
three variables over
$\FF_q$ cutting out a possibly singular plane model of $C$ within
the projective plane $\PP^2$, i.e., a projective, geometrically
irreducible one-dimensional scheme 
$C'$ whose normalization is isomorphic to $C$.
Let $d$ be the degree of the polynomial; then by 
Pl\"ucker's adjunction
formula, the genus $g$ of $C$ is at most $(d-1)(d-2)/2$. That is,
$g$ is bounded by a polynomial in $d$. We will assume also conversely that
$d$ is bounded by a polynomial in $g$, so that polynomiality can be measured in 
terms of $d$ rather than $g$. This is no real restriction: by
Riemann-Roch, any curve of degree $g$ admits a singular plane model of degree $g$, so can be properly input into our algorithm.

We need to explicitly describe the singularities of 
$C'$ and the sequence of blowups of $\PP^2$ that resolves these
singularities. Straightforward algorithms for doing this require passing
to extensions of $\FF_q$ whose degree is not polynomial in the input length
(e.g., an extension over which all singular points become rational).
However, there exist methods that perform the resolution of singularities
in polynomial time, e.g., that of Kozen \cite{kozen}. Note that the number
of $\overline{\FF_q}$-rational points of $C$ lying over singular points
of $C'$ is at most $(d-1)(d-2)/2$, since each one contributes at least one
to the discrepancy between the Pl\"ucker bound and $g$.

Put $m = \lceil 2 \log_q(d) \rceil$.
Since there are at most $(d-1)(d-2)/2$ geometric points of $C$
lying above singular points on $C'$, we can draw an $\FF_{q^m}$-rational
line in $\PP^2$ not meeting any of the singular points. Pick such a line, let $F$ be the divisor in which
the line meets $C'$, and choose an $\overline{\FF_q}$-point $O$ of $F$; then
$O$ is defined over $\FF_{q^{mn}}$ for some $n \leq d$.

The key to constructing a black box presentation of $\Cl(C)$ is that
the Riemann-Roch theorem on $C_{mn}$
can be made (Las Vegas) polynomial time effective; in other words,
given a divisor $D$ on $C_{mn}$, 
one can efficiently test functions
for membership in $L(D)$,
write down a basis of $L(D)$,  and express elements of $L(D)$ as linear
combinations of that basis.
See for instance 
Huang and Ierardi \cite[\S 2]{huang-ierardi} for an explicit construction;
see also Volcheck \cite{volcheck1}, \cite{volcheck2} for a somewhat more
practical construction.
(Note that Huang and Ierardi assume that all singular points are rational,
but they also point out that this restriction is only needed to ensure
that resolution of singularities can be performed efficiently. Thanks
to the argument of Kozen from \cite{kozen}, this restriction can be lifted.)

Let $S$ be the set of effective
divisors $E$ on $C_{mn}$
with $\deg(E) \leq g$ and $\ell(E) = 1$,
represented as bit strings by listing the $\overline{\FF_q}$-points
on $E$ (on the blowup of $\PP^2$ chosen to resolve the singularities
of $C'$). Then given a divisor $D$ of degree 0, we can describe
a reduction   procedure to produce $E \in S$ with $D \sim
 E - \deg(E)(O)$ as follows. Apply effective Riemann-Roch to produce
$E_0$ of degree $g$ with $E_0 \sim D + g(O)$. Then
repeatedly apply 
effective Riemann-Roch to find divisors $E_1, E_2, \dots$ with
$\deg(E_i) = g - i$ and $E_i - (g-i)(O) \sim E_{i+1} - (g-i-1)(O)$,
until it is no longer possible to do so. 
If this stops at $E_i$, then $E_i \in S$
and $E_i - (g-i)(O) \sim D$.

To add $D_1, D_2 \in S$,
we may apply the reduction procedure to $D_1 + D_2 - \deg(D_1+D_2)(O)$.
To negate $D \in S$, we may apply the reduction procedure to
$-D + \deg(D)(O)$.
Hence we have produced a black box presentation with unique
encodings $f: S \to \Cl(C_{mn})$, modulo the problem of exhibiting
a generating set. We discuss generating sets in the next section.

\section{Finding generators of class groups}

With notation as in the previous section,
let $T$ be the subset of $S$
corresponding to elements of $\Cl(C)$. In order to have a black box presentation with unique encodings $f: T \to \Cl(C)$, 
so that we can apply Watrous's algorithm to compute $\#\Cl(C)$, we need to
exhibit with high probability a subset of $T$ which generates $\Cl(C)$; to do this
provably (without too much headache), we will have to assume that $q$ is ``not too small''
compared to $g$. It may be possible to lift this restriction with an even
more elaborate argument than the already involved procedure given below.

We first observe that it suffices to somehow generate uniformly random
 elements of 
$\Cl(C)$.
\begin{lemma} \label{L:random}
Let $G$ be a finite abelian group of order $\leq 2^h$. Then for any nonnegative integer $i$, if one chooses $h+i$ elements of $G$ uniformly at random (with 
replacement), the probability that the chosen elements generate $G$ is at least
$1 - 2^{-i}$.
\end{lemma}
\begin{proof}
As stated, this is \cite[Theorem~D.1]{lomont}; the argument therein is due
to Pak \cite{pak1}, \cite{pak2}. (Roughly, one checks that the probability is minimized by elementary 2-groups, then verifies the bound explicitly in that case.) An older but weaker result in the same spirit (which only yields the desired probability $1-2^{-i}$ after sampling on the order of $2h+i$ elements, rather than $m+i$) is due to Erd\H{o}s and R\'enyi
\cite[Theorem~1]{erdos-renyi}.
\end{proof}

By a \emph{$b$-uniform oracle} on a finite set $V$,
we will mean an oracle which either fails to return an answer with probability at most
$1/4$, or returns a element of $V$ according to a 
probability distribution $p: S \to [0,1]$ such that
for any $x,y \in V$, $p(x) \leq bp(y)$. (The constant $1/4$ is chosen merely for
definiteness; as in Section~\ref{sec:prob}, there is no harm
in replacing $1/4$ by any other fixed constant between 0 and 1.)

\begin{lemma} \label{L:uniform div}
Given a positive integer $e$ such that $16g < q^{e/2}$,
let $V$ be the set of prime divisors on $C$ of degree $e$. Then 
there exists a $(1 + (2g-2+d)/e)$-uniform
oracle on $V$, running in time polynomial in $g$ and $\log(q)$.
\end{lemma}
\begin{proof}
Put $j= \lceil(2g-1+e)/d \rceil$.
Consider an oracle that performs the following operation: select a random
homogeneous polynomial over $\FF_q$ of degree $j$, then
extract uniformly at random an $\FF_{q^e}$-rational point of $C$
on which this polynomial vanishes, and return the divisor consisting of the
Galois orbit of that point. (Here failure occur if there is no such point, if the chosen polynomial restricts to zero on $C$, or if Las Vegas univariate
polynomial factorization fails.)

To analyze this oracle, we first note that the homogeneous polynomials
of degree $j$ give rise to $q^{jd+1-g}$ distinct functions on $C$, by
Riemann-Roch (and each occurs the same number of times). Also by Riemann-Roch,
each prime divisor $E$ of degree $e$ occurs in the zero locus of $q^{jd-e+1-g}$
such functions: namely, if $F$ is
the divisor along which $C$ meets some line, we have
$\ell(jF - E) = jd-e + 1-g + \ell(K-jF+E) = jd-e+1-g$ since
$\deg(K-jF+E)=2g-2-jd+e < 0$.

Note that each nonzero homogeneous polynomial of degree $j$ can give rise
to at most $\lfloor jd/e \rfloor$ distinct divisors. This means that on
one hand, the ratio between the probabilities of producing any two prime divisors of 
degree $e$ is at most $\lfloor jd/e \rfloor \leq 1 + (2g-2+d)/e$.
On the other hand, by Lemma~\ref{L:count primes}
 and the hypothesis $16g < q^{e/2}$,
the probability of success of the oracle (assuming success in the polynomial factorization,
which can be assured to sufficiently high probability by repeated trials) is 
at least
\begin{align*}
\frac{1}{e} \frac{(q^{jd+1-g-e}-1) (q^e(1-q^{-1}) - 4gq^{e/2})}{q^{jd+1-g}} 
&> \frac{1}{2e} \frac{q^e (1-q^{-1}) - 4g q^{e/2}}{q^e} \\
&\geq \frac{1}{4e} \frac{q^e - 8g q^{e/2}}{q^e} \\
&> \frac{1}{8e}.
\end{align*}
With $1024e^2$ invocations of this oracle (as in Section~\ref{sec:prob}),
we can boost this probability to $3/4$,
yielding the desired result.
\end{proof}
Note that one cannot state the previous lemma as written without some lower bound
on $q$ with respect to $g$; otherwise it might happen that $V$ is empty, and one
certainly cannot construct the desired oracle in that case! This complication is
the reason we will be limited to the case where $q$ is ``not too small'' below.

Next, we give a ``simulation'' conversion from a 
$b$-uniform oracle into a $1$-uniform oracle. (The ``simulation'' qualifier refers
to the fact that one must explicitly know the probability distribution on the initial
oracle, which is too strong an assumption to make in practice.)
\begin{lemma} \label{L:convert}
Suppose we are given a $b$-uniform oracle on a finite set $V$ with \emph{known} distribution
and error probability. Then we can construct a $1$-uniform oracle on $V$ requiring
at most $16b^2$ invocations of the initial oracle.
\end{lemma}
\begin{proof}
Let $p: V \to [0,1]$ be the probability distribution of the initial oracle,
and put $p_0 = \min_{x \in V} \{p(x)\}$; note that
\[
p_0 \geq \frac{1}{b \#V}
\]
since the initial oracle is $b$-uniform.
Consider the following operation: invoke the initial oracle once to produce $x$,
then return $x$ with probability $p_0/p(x)$ and fail otherwise. This operation
is equally likely to return any element of $V$, and succeeds with probability
$p_0 \#V \geq 1/b$. Performing the operation $16b^2$  times (as in Section~\ref{sec:prob})
gives a new oracle with failure probability at most $1/4$, as desired.
\end{proof}

We now put together the previous lemmas.
It should be cautioned that the awkward intricacy of the resulting
Lemma~\ref{L:generators} is
caused by our desire to have a fully unconditional
 complexity analysis; in practice,
one is quite likely to obtain a generating set by selecting divisors by any
reasonably arbitrary process!

\begin{lemma} \label{L:generators}
Under the assumption $16g < q^{1/2}$,
there exists a Monte Carlo algorithm that
produces a subset of $T$ generating $\Cl(C)$ in time polynomial in $g, \log(q)$.
\end{lemma}
\begin{proof}
Put $h = \lceil \log_2 (q^{g/2} (\sqrt{q}+1)^g) \rceil$, so that
$\#\Cl(C) \leq 2^h$ by Proposition~\ref{P:Weil bound}.
Put $N = \lceil (1 + (2g-2+d)/e) \rceil$.
By repeated use of Lemma~\ref{L:uniform div} together with
Lemma~\ref{L:convert}, we can produce, for each of
$i=1, \dots, 2g+1$,
a list of $32N^2 (2g-1)^2(h+3)$ prime divisors of degree $i$, each produced
by an $N$-uniform oracle.
Moreover, we can do this with overall probability of failure at most 1/16.

Apply Lemma~\ref{L:uniform div} to produce a divisor $U$ of degree 1,
then convert each divisor $E$ in the list into an element of $T$ by reducing
$E - \deg(E) U$ via effective Riemann-Roch. We now verify that the resulting
elements of $\Cl(C)$ generate $\Cl(C)$ with probability at least $3/4$
by a ``simulation'' argument; namely, we exhibit the existence of another random
process which necessarily produces a sublist of our given list, but which
also produces a generating set for $\Cl(C)$ with probability at least 3/4.

In this context, we may assume that we know
the distribution of the oracle produced by Lemma~\ref{L:uniform div}. (In the context
of constructing the algorithm, we cannot use this knowledge, as it would amount to
already knowing  the zeta function of $C$. The point is 
that we do not use the information to perform any algorithmic steps, only to verify
the error bound.) By Lemma~\ref{L:convert}, we may then extract from the given data
a list of $(2g-1)(h+3)$ prime divisors of degree $i$, with failure probability
at most $1/16$. (The factor of $16N^2$ is shed in the application of Lemma~\ref{L:convert};
shedding the factor of $2(2g-1)$ allows us to shrink the failure probability to
$1/16^{2g-1}$, so that the combined failure probability after producing all $2g-1$ lists is
at most $1/16$.)

From these new lists, we can in turn simulate the uniform random 
choice of $h+3$ divisors of degree $2g-1$. We do this assuming knowledge of the
number of prime divisors of degrees $1,\dots,2g-1$ (again, this amounts to knowing
the desired zeta function, but this is okay for proving an error bound). With that
knowledge, we may choose a ``shape'' of a degree $2g-1$ divisor (i.e., the information
of how many prime divisors occur with a given degree and multiplicity) 
according to the distribution
which is uniform for individual divisors.
(That is, each shape has probability proportional
to the number of divisors taking that shape.)
Given a shape, we may then read off from
our lists uniformly random prime divisors of the appropriate lengths; we cannot use more than
$2g-1$ divisors of any one length at a time, so we have enough data to do this
$h+3$ times (with no additional failure probability at this step).

Finally, with $h+3$ uniformly random divisors of degree $2g-1$ in hand, we obtain by
reduction $h+3$ uniformly random elements of $\Cl(C)$ (by the calculations
of Section~\ref{sec:rep}). By  Lemma~\ref{L:random}, these generate $\Cl(C)$
with probability at least $1 - 1/8$. Since the divisors we produced were synthesized
from the original list we produced, that list also generates $\Cl(C)$ with probability
at least $1-1/8$. Totaling the failure and error probabilities yields an
error probability in the Monte Carlo algorithm of $1/4$, as desired. (Note that the
only step which is Monte Carlo rather than Las Vegas is the last one, since we do not check whether the random
elements we produced actually do generate $\Cl(C)$.)
\end{proof}

We now may combine all of our efforts so far to obtain the following result.
\begin{prop} \label{P:quantum alg}
For $e$ such that $16g < q^{e/2}$, there
exists a quantum algorithm to compute
$\#\Cl(C_e)$ in time polynomial in $g,\log(q),e$.
\end{prop}
\begin{proof}
The construction of the previous section exhibits a black box presentation with
unique encodings for $\Cl(C_e)$, minus the construction of a set of generators;
these are furnished by Lemma~\ref{L:generators}. Now Lemma~\ref{L:quantum}
applies to yield the desired algorithm.
\end{proof}

\section{Computing the zeta function}
\label{sec:zeta}

Retain notation as in Section~\ref{sec:Jacobians}.
By Proposition~\ref{P:quantum alg}, we can exhibit a quantum algorithm
to compute the order of the group $\#\Cl(C_n)$ in time
polynomial in $g, \log(q), n$, as long as $16g < q^{n/2}$.
With this quantum input in hand, we now establish Theorem~1.
\begin{proof}[Proof of Theorem~1]
We first proceed under the assumption that $16g < q^{1/2}$, so that
we may apply Proposition~\ref{P:quantum alg} for any $e$. Note that
this assumption only intervenes via the invocation of 
Proposition~\ref{P:quantum alg}; if one were to prove a form of that
proposition without the lower bound on $q$, this restriction would
drop out of the proof of Theorem~1.

Recall
that by the Weil conjectures (see Proposition~\ref{P:Weil bound} and 
also \cite[Appendix~C]{hartshorne}),
we can factor $P(t)$ over $\CC$ as
\[
(1 - r_1 t)\cdots(1 - r_{2g} t),
\]
where each $r_i$ is an algebraic integer of absolute value $q^{1/2}$,
and $r_i r_{g+i} = q$ for $i=1, \dots, g$. 
Write $P(t) = a_0 + a_1 t + \cdots + a_{2g} t^{2g}$ with $a_0 = 1$; then 
the symmetry $r_i r_{g+i} = q$ implies that
$a_{g+i} = q^i a_{g-i}$ for $i=1, \dots, 2g$, so to determine $P(t)$
it is enough to determine the integers $a_1, \dots, a_g$.

As noted earlier (Proposition~\ref{P:Jacobian}), we then have
\[
\#\Cl(C_n) = \prod_{i=1}^{2g} (1 - r_i^n) = q^{gn} \prod_{i=1}^{2g}
(1 - r_i^{-n}).
\]
Put
\[
c_n = q^{-gn} \#\Cl(C_n), \qquad s_n = q^{-n} \frac{1}{n} \sum_{i=1}^{2g}
r_i^n = \frac{1}{n} 
\sum_{i=1}^{2g} r_i^{-n};
\]
then we can write
\[
-\frac{\log c_n}{n} = \sum_{j=1}^\infty s_{nj}.
\]
By the Newton-Girard formulae,
\[
n q^n s_n + a_1 (n-1) q^{n-1} s_{n-1} + \cdots + a_{n-1} q s_1 + n a_n = 0
\qquad 
(n=1, \dots, g);
\]
in particular, it is enough to determine $s_1, \dots, s_g$, as we can then
recover $a_1, \dots, a_g$.

Using Proposition~\ref{P:quantum alg}, 
we can compute $c_n$ in suitable time for $n=1, \dots, m$ with $m = \max\{18,2g\}$.
We can then compute $s_1, \dots, s_g$ exactly as follows.
Suppose $n \leq g$ and that
$s_i$ has been computed exactly for $i=1, \dots, n-1$.
By the Newton-Girard formulae, the residue modulo $n$ of the integer
$n q^n s_n$ is determined by
$s_1, \dots, s_{n-1}$.
Hence we can recover the exact value of $s_n$ if we can compute $q^n s_n$
to within an error of less than $0.5$.

Let $\mu(n)$ denote the M\"obius function, put $k = \lfloor m/n \rfloor$,
and compute
\[
q^n e_n = \sum_{i=1}^{k}
-q^n \mu(i) \frac{\log c_{ni}}{ni}
= q^n s_n + \sum_{j=k+1}^\infty q^n a_{n,j} s_{nj}
\]
to an error of less than $0.005$. Here
\[
a_{n,j} = \sum_{1 \leq i \leq k, i | j} \mu(i)   
\] 
is an integer of absolute value at most $k$,
so
\begin{align*}
\left| \sum_{j > k} q^n a_{n,j} s_{nj} \right|
&\leq q^n \sum_{j=k+1}^\infty \frac{2gk q^{-nj/2}}{nj} \\
&\leq q^n \frac{2gk}{n} \sum_{j=k+1}^\infty \frac{q^{-nj/2}}{k+1} \\
&= \frac{2g}{n} \frac{q^{-n(k-1)/2}}{1-q^{-n/2}} \\
&\leq (k+1) \frac{q^{-n(k-1)/2}}{1-q^{-n/2}}.
\end{align*}
This last expression is less than $0.495$ if $k \geq k_0$ and $q^n \geq q_0$
for each of
\[
(k_0, q_0) \in \{ (2,50), (3,14), (4,7), (5,5), (6,4), (8,3), (15,2)\}.
\]
Note that $18 \geq (k_0+1)\log_2(q_0)$ for each pair $(k_0,q_0)$ in the
above list. Since $m \geq 18$ and $q \geq 2$, for any pair $(k,n)$ with
$k \geq 2$ and $k = \lfloor m/n \rfloor$, we then
have $k \geq k_0$ and $q^n \geq q_0$ for some pair $(k_0,q_0)$.
Thus the computed value of $q^n e_n$ differs from
$q^n s_n$ by less than $0.5$, so we may determine $s_n$ exactly.
We may thus recover the zeta function in this fashion.

To recap, we have proved that we can recover the zeta function of $C$ provided
that $16g < q^{1/2}$; it remains to relax this restriction.
Given arbitrary $g$ and $q$, choose $m_1, m_2$ subject to
the following conditions.
\begin{itemize}
\item $m_1 < m_2$.
\item For $i=1,2$, $m_i$ is prime and $m_i-1$ is divisible by some prime greater than
$2g$.
\item $16g < q^{m_1/2}$.
\end{itemize}
The existence of such $m_1, m_2$ of size bounded by a polynomial in $g, \log(q)$
is guaranteed, e.g., by a theorem of Harman 
\cite[Theorem~1.2]{harman}, which asserts that for any fixed
$\theta \leq 0.610$, there exist effectively computable constants $\delta>0$
and $x_0 \in \mathbb{R}$ such that for $x \geq x_0$,
there are at least $\delta x/\log(x)$ primes $p \in \{1,\dots, x\}$ such that
$p-1$ has greatest
prime factor bigger than $x^\theta$. (Many  results of this ilk
exist in the analytic number theory
literature,
but the effective computability of the constants seems to be new to \cite{harman}.)

Apply the previous argument to compute the zeta functions of $C_{m_1}, C_{m_2}$.
We thus have the lists $r_1^{m_1}, \dots, r_{2g}^{m_1}$
and $r_1^{m_2},\dots, r_{2g}^{m_2}$. By the construction of $m_1$ and $m_2$,
the field extension $\QQ(r_1,\dots, r_{2g})$ cannot contain a nontrivial $m_1$-st
or $m_2$-nd root of unity (else such a root of unity would generate a field whose degree
contains a prime factor greater than $2g$, whereas the degree of
$\QQ(r_1, \dots, r_{2g})$ divides $(2g)!$). Thus we have
$(r_j^{m_1})^{m_2} = (r_l^{m_1})^{m_2}$ if and only if $r_j = r_l$.

If we now pick out an element $A$ of the first list, there is only
one value (possibly repeated) $B$ occurring in the second list with
$A^{m_2} = B^{m_1}$. We can thus unambiguously (up to interchanging identical values)
pair off each $r_j^{m_1}$ with its corresponding $r_j^{m_2}$, and then recover the
$r_j$. This completes the proof.
\end{proof}

\section{Cyclic resultants}

The above argument can also be described as follows. Given a polynomial
$P(t)$ with roots $r_1, \dots, r_d$, the \emph{$m$-th cyclic resultant} 
of $P(t)$ is defined as
\[
\Res(P(t), t^m - 1) = \prod_{i=1}^d (r_1^m - 1).
\]
These arise in a number of applications; see \cite{hillar} for further
discussion. A theorem of Fried \cite{fried} asserts that if $P(t)$
has even degree and is reciprocal (i.e., $P(t) = t^d P(1/t)$),
then $P$ is uniquely determined by its sequence of cyclic
resultants. This is precisely the situation in which we are in, which is not
surprising: Fried arrived at this situation by counting fixed points of the
powers of an endomorphism of a topological torus in terms of the Lefschetz
trace formula on cohomology, and we are doing the same with the Frobenius
endomorphism on an abelian variety.

Unfortunately, Fried's theorem does not give an effective bound on the
number of cyclic resultants needed to recover $P(t)$, nor an algorithm for
doing so. A conjecture of Sturmfels and Zworski asserts that
the first $d/2 + 1$ cyclic resultants should suffice for $P$ generic
(if $P$ is not reciprocal, they conjecture that generically $d+1$ resultants suffice).
A theorem of Hillar and Levine \cite{hillar-levine} states that the first
$2^{d+1}$ cyclic resultants determine $P$; what we have done is show that
for very special reciprocal $P$, we can explicitly recover $P$ from only
$d$ cyclic resultants.

Whether one can bring $d$ down any closer to the theoretical lower bound
$d/2$, i.e., whether one can compute the zeta function of a curve of genus
$g$ using fewer than $2g$ calls to the quantum oracle, is a tantalizing
question. Our current approach fails to accomplish this because, for 
instance, we recover $s_g$ from $s_g + s_{2g} + \cdots$, and the term $s_{2g}$
is of exactly the same order as the size of the interval in which
we must bound $s_g$ in order to determine it exactly, namely $q^{-g}$. 
Thus breaking the
$2g$ barrier would seem to require a fundamental new idea.

Incidentally,
this barrier may be of interest even in the absence of quantum computers,
as it may be possible to use the proof
of Theorem~1 to obtain a probabilistic polynomial time algorithm
for verifying the zeta function of a curve, which verifies the orders
of the first few Jacobian groups. Unfortunately, while it is easy to
efficiently verify the exponent of a black box group, it is less clear how
to efficiently verify its order. (Thanks to Dan Bernstein for this remark.)

\section{Further comments}

It should be noted that the problem of giving an efficient quantum algorithm to
compute the zeta function of an arbitrary variety $X$ over a finite field
$\FF_q$ is now effectively solved in dimension $\leq 1$. For $\dim (X) = 0$,
i.e., for $X$  a finite union of closed points, computing the zeta function
of $X$ amounts to finding the distinct-degree factorization of a monovariate
polynomial, so this can even be done in deterministic polynomial time. For
$\dim(X) = 1$, if $X$ is geometrically irreducible, one can find the unique
smooth projective curve $C$ birational to $X$, compute its zeta function,
then express the discrepancy between the zeta functions of $X$ and $C$
in terms of the zeta functions of zero-dimensional varieties. If $X$
is not geometrically irreducible, one can split it over an extension of
degree at most its genus and proceed as above.

However, considering varieties of a fixed higher dimension seems to pose
more serious challenges. (Allowing the dimension to vary brings us dangerously
close to the $P = NP$ problem, which we prefer to stay well clear of.)
Things are well understood, at least theoretically, if the characteristic
$p$ of $\FF_q$ is fixed; as noted earlier, Lauder and Wan \cite{lauder-wan}
give a deterministic
algorithm for computing the zeta function of a singular hypersurface
of degree $d$ in $\PP^n$, in time polynomial in $p, \log_p(q), d$.
(Again, one can reduce to this case by induction on dimension, since
any irreducible variety is birational to a hypersurface.)

On the other hand, if $p$ is allowed to vary, then even the following question
remains somewhat mysterious, except in some cases related to modular forms
(as demonstrated by ongoing work of Bas Edixhoven and his collaborators on
efficient computation of the values of Ramanujan's $\tau$ function).
\begin{question}
Let $X$ be a fixed variety over $\QQ$ (or better, fix a model over $\ZZ$)
of dimension greater than $1$. 
Does there necessarily exist a deterministic,
random, or quantum polynomial time algorithm in $\log(p)$ to determine
the zeta function of $X$ over $\FF_p$, for $p$ a varying prime?
\end{question}
For $X$ of dimension 1, Schoof-Pila gives a deterministic affirmative answer.
However, the approach used there breaks down in higher dimensions; briefly
put, there is no ``geometric'' realization of the higher \'etale cohomology groups
analogous to the realization of the first \'etale cohomology group in the
Tate module of the Jacobian. The work of Edixhoven suggests such a realization
in case the relevant cohomology group is ``modular'', by comparing the
higher \'etale cohomologies to first \'etale cohomologies on other spaces.
However, already the case when $X$ is a (fixed) surface of
general type, without any special structure, seems to require a new idea.

We also point out a related  but markedly different investigation initiated by
van Dam \cite{vandam}, who looks for ``efficient'' quantum circuits 
for computing the zeta functions of varieties, mostly in dimensions
greater than 1. The emphasis there is on directly realizing Frobenius
eigenvalues within easy-to-construct Hermitian operators; this is done in
\cite{vandam} for some diagonal hypersurfaces (where the relevant eigenvalues
are Gauss sums) but seems quite difficult in general.

\subsection*{Acknowledgments}

Thanks to Sean Hallgren for providing helpful
comments on a prior version of the manuscript,
to Igor Pak for providing the reference \cite{lomont},
and to Laci Babai for some clarifications about black box groups.
The author is funded by NSF grant DMS-0400747.

\end{document}